\documentclass{amsart}
\usepackage{amssymb}
\newtheorem{thm}{Theorem}
\newtheorem{lem}{Lemma}
\newtheorem{prop}{Proposition}

\newtheorem{df}{Definition}

\begin{document}

\bibliographystyle{plain}

\title[ Romi F. Shamoyan]{On multipliers of mixed -norm holomorphic $F^{p,q}_\alpha$ type spaces on the unit polydisc}

\author[]{Romi F. Shamoyan}
\author []{Milos Arsenovic}

\address{Bryansk University, Bryansk Russia}
\email{\rm rshamoyan@gmail.com}
\address{Belgrade University,Serbia}
\email{\rm arcenovis@matf.bg.ac.rs}

\date{}

\maketitle

\begin{abstract}
We describe certain new spaces  of coefficient multipliers of analytic Lizorkin-Triebel $F^{p,q}_\alpha$ type spaces in the unit polydisc with some restriction on parameters. This extends some previously known assertions on coefficient multipliers in classical Bergman $A^{p}_{\alpha}$ spaces in the unit disk.

\end{abstract}

\footnotetext[1]{Mathematics Subject Classification 2010 Primary 30H10.  Key words
and Phrases: Multipliers, polydisc, analytic functions, analytic Lizorkin-Triebel spaces.}

\section{Introduction and preliminary results}

The goal of this paper is to continue the investigation of spaces
of coefficient multipliers of analytic $F^{p,q}_\alpha$
Lizorkin-Triebel type spaces in the unit polydisc started  in
\cite {Sh1}, \cite{Sh2}. These $F^{p,q}_\alpha$ type spaces serve
as a very natural extension of the classical Hardy and Bergman
spaces in the unit polydisc simultaneously. Spaces of Hardy and
Bergman type in higher dimension were studied intensively by many
authors during past several decades, see for example
\cite{Sv1},\cite{DS},\cite{GL},\cite {OF1},\cite{OF2},\cite{OF3}
and references therein. The investigation in this new direction of
$F^{p,q}_\alpha$ type classes in the unit polydisc was started in
recent papers of the second author, see \cite{Sh1}, \cite{Sh2},
\cite{Sh3}. We note that complete analogues of these classes in
the unit ball were studied also in recent papers of Ortega and
Fabrega see \cite{OF1}, \cite{OF2}, \cite{OF3} and references
therein. See also \cite{SSV}, \cite{GL}, \cite{SL} for some other
properties of these new analytic $F^{p,q}_\alpha$ type classes and
related spaces in the unit disk or subframe. Below we list
notation and definitions which are needed and in the next section
we formulate and prove main results of this note.

Let $\mathbb D = \{ z \in \mathbb C : |z| < 1 \}$ be the unit disc in $\mathbb C$, $\mathbb T = \partial \mathbb D =
\{ z \in \mathbb C : |z| = 1 \}$, $\mathbb D^n$ is the unit polydisc in $\mathbb C^n$, $\mathbb T^n \subset \partial
\mathbb D^n$ is the distinguished boundary of $\mathbb D^n$, $\mathbb Z_+ = \{ n \in \mathbb Z : n \geq 0 \}$, $\mathbb Z_+^n$ is the set of all multi indexes and $I = [0, 1)$.

We use the following notation: for $z = (z_1, \ldots, z_n) \in \mathbb C^n$ and $k = (k_1, \ldots, k_n) \in
\mathbb Z_+^n$ we set $z^k = z_1^{k_1} \cdots z_n^{k_n}$ and for $z = (z_1, \ldots, z_n) \in \mathbb D^n$ and $\gamma \in
\mathbb R$ we set $(1-|z|)^\gamma = (1-|z_1|)^\gamma \cdots (1-|z_n|)^\gamma$ and $(1-z)^\gamma = (1-z_1)^\gamma \cdots
(1-z_n)^\gamma$. Next, for $z \in \mathbb R^n$ and $w \in \mathbb C^n$ we set $wz = (w_1z_1, \ldots, w_nz_n)$. Also, for $k \in \mathbb Z_+^n$ and $a \in \mathbb R$ we set $k + a = (k_1 + a, \ldots, k_n + a)$. For $z = (z_1, \ldots, z_n) \in \mathbb C^n$ we set $\overline z = (\overline z_1, \ldots, \overline z_n)$. For $k \in (0, +\infty)^n$ we set $\Gamma (k) = \Gamma (k_1) \cdots \Gamma (k_n)$.

The Lebesgue measure on $\mathbb C^n \cong \mathbb R^{2n}$ is denoted by $dV(z)$, normalized Lebesgue measure on $\mathbb T^n$
is denoted by $d\xi = d\xi_1 \ldots d\xi_n$ and $dR = dR_1 \ldots dR_n$ is the Lebesgue measure on $[0, +\infty)^n$.

The space of all functions holomorphic in $\mathbb D^n$ is denoted by $H(\mathbb D^n)$. Every $f \in H(\mathbb D^n)$ admits
an expansion $f(z) = \sum_{k \in \mathbb Z_+^n} a_k z^k$. For $\beta > -1$ the operator of fractional differentiation is
defined by
\begin{equation}\label{DfFrD}
D^\beta f(z) = \sum_{k \in \mathbb Z_+^n} \frac{\Gamma(k+\beta+1)}{\Gamma(\beta+1)\Gamma(k+1)} a_k z^k,
\qquad z \in \mathbb D^n.
\end{equation}
For $f \in H(\mathbb D^n)$, $0 < p < \infty$ and $r \in I^n$ we set
\begin{equation}\label{EqMpr}
M_p(f,r) = \left( \int_{\mathbb T^n} |f(r\xi)|^p d\xi \right)^{1/p},
\end{equation}
with the usual modification to include the case $p = \infty$. For $0<p\leq \infty$ we have analytic Hardy classes in the
polydisc:
\begin{equation}\label{EqDfHp}
H^p(\mathbb D^{n}) = \{ f \in H(\mathbb D^n) : \| f \|_{H^p} = \sup_{r \in I^n} M_p(f,r) < \infty \}.
\end{equation}

For $n=1$ these spaces are well studied. The topic of multipliers of Hardy spaces in polydisc is relatively new, see for example \cite{Sv1}, \cite{Tr1}. These spaces are Banach spaces for $1\leq p$ and complete metric spaces for all other positive values of $p$. Also, for $0<p\leq \infty$, $0<q<\infty$ and $\alpha > 0$ we have mixed (quasi) norm spaces, defined below.

\begin{equation}\label{EqMixN}
A^{p,q}_\alpha (\mathbb D^n) = \left\{ f \in H(\mathbb D^n) : \| f \|^q_{A^{p,q}_\alpha} = \int_{\mathbb I^n}
M_p^q(f, R) (1-R)^{\alpha q - 1} dR < \infty \right\}.
\end{equation}

These spaces are Banach spaces for cases when $1\leq\min{p,q}$, and they are complete metric spaces for all other values of parameters. We refer the reader for these classes in the unit ball and the unit disk to \cite{Sv1}, \cite{OF1}, \cite{OF2}, \cite{OF3} and references therein. Multipliers between $A^{p,q}_\alpha$ spaces on the unit disc were studied in
detail in \cite{JP}.

As is customary, we denote positive constants by $C$, sometimes we indicate dependence of a constant on a parameter by
using a subscript, for example $C_q$.

We define now the main objects of this paper.

For $0<p,q<\infty$ and $\alpha > 0$ we consider Lizorkin-Triebel spaces $F^{p,q}_\alpha (\mathbb D^n) = F^{p,q}_\alpha$
consisting of all $f \in H(\mathbb D^n)$ such that
\begin{equation}\label{EqTrLi}
\| f \|^p_{F^{p,q}_\alpha} = \int_{\mathbb T^n} \left( \int_{I^n} |f(R\xi)|^q (1-R)^{\alpha q - 1}dR \right)^{p/q} d\xi < \infty.
\end{equation}

It is not difficult to check that those spaces are complete metric spaces, if $\min(p,q) \geq 1$ they are Banach spaces.

We note that this scale of spaces includes, for $p=q$, weighted
Bergman spaces $A^p_\alpha = F^{p,p}_{\frac{\alpha+1}{q}}$, see
\cite{DS} for a detailed account of these spaces. On the other
hand, for $q=2$ these spaces coincide with Hardy-Sobolev spaces
namely $H^p_\alpha = F^{p,2}_{\frac{\alpha + 1}{2}}$, for this
well known fact see \cite {Sv1}, \cite{SSV} and references
therein.

Finally, for $\alpha \geq 0$ and $\beta \geq 0$ we set
\begin{equation}\label{EqDinin}
A^{\infty,\infty}_{\alpha, \beta} (\mathbb D^n) = \{ f \in H(\mathbb D^n) : \| f \|_{A^{\infty,\infty}_{\alpha,\beta}} =
\sup_{r \in I^n} (M_\infty (D^\alpha f, r)) (1-r)^\beta < \infty \}.
\end{equation}

This space is a Banach space. For all positive values of $p$ and $s$ we introduce the following two new spaces.

The limit space case $F^{p,\infty,s}(\mathbb D^n)$ is defined  as a space of all analytic functions $f$ in the polydisc
such that the function $\phi(\xi) = \sup_{r \in I^n}|f(r\xi)|(1-r)^{s}$, $\xi \in \mathbb T^n$ is in $L^p(\mathbb T^n, d\xi)$.

Finally, the limit space case $A^{p,\infty,s}(\mathbb D^n)$ is the space of all analytic functions $f$ in the polydisc such that
$$\sup_{r\in I^n}M_{p}(f,r)(1-r)^{s}<\infty .$$

Obviously the limit $F^{p,\infty,s}$ space is embedded in the last space we defined. This simple observation will be used at the end of this note. These last two spaces are Banach spaces for all $p\geq 1$ and they are complete metric spaces for all other positive values of $p$.

The following definition of coefficient multipliers is well known in the unit disk. We provide a natural extension to the
polydisc setup.
\begin{df}
Let $X$ and $Y$ be  quasi normed subspaces of $H(\mathbb D^n)$. A sequence $c = \{ c_k \}_{k \in \mathbb Z_+^n}$ is said to be a coefficient multiplier from $X$ to $Y$ if for any function $f(z) = \sum_{k \in \mathbb Z_+^n} a_k z^k$ in $X$ the function $h = M_cf$ defined by $h(z) = \sum_{k \in \mathbb Z_+^n} c_k a_k z^k$ is in $Y$. The set of all multipliers from $X$ to $Y$ is denoted by $M_T(X, Y)$.
\end{df}

The problem of characterizing the space of multipliers (pointwise multipliers and coefficient multipliers)
between various spaces of analytic functions is a classical problem in complex function theory,
there is vast literature on this subject, see \cite{BP}, \cite{DS}, \cite{JP}, \cite{MP1}, and references therein.

In this paper we are looking for some extensions of already known theorems, namely we are interested in spaces of multipliers acting into analytic Lizorkin-Triebel $F^{p,q}_\alpha$ type spaces in the unit polydisc and from these spaces into certain well studied classes like mixed norm spaces, Bergman spaces and Hardy spaces. We note that the analogue of this problem of description of spaces of multipliers in $\mathbb R^n$ for various functional spaces in $\mathbb R^n$ was considered previously by various authors in recent decades, for these we refer the reader to \cite{Tr2}.

\section{On coefficient multipliers of analytic Lizorkin - Triebel type spaces in polydisc}

We start this section with several lemmas which will play an important role in the proofs of all our results. We note that these assertions can serve as direct natural extensions of previously known one dimensional results to the case of several complex variables. The following lemma in dimension one is well-known, see \cite{Sh1}, \cite{DS}, \cite{Ga1} and references therein. This is a several variables version of the so called Littlewood-Paley formula, see \cite{Ga1}, which allows to pass from integration on the unit circle to integration over unit disk. We omit a simple proof which follows from a straightforward technical computation based on orthonormality of the trigonometric system, similarly to the classical planar case.

\begin{lem}[\cite{Sh2}]\label{LeCon}
For $f, g \in H(\mathbb D^n)$, $r \in I^n$ and $\alpha > 0$ we have
\begin{align}\label{EqCon}
\int_{\mathbb T^n} f(r\overline t) g(rt) dt = & (2 \alpha) ^n \prod_{j=1}^n r_j^{-2\alpha} \times  \int_0^{r_1}  \cdots \int_0^{r_n} \cr
& \int_{\mathbb T^n} D^{\alpha + 1} g(R\xi) f(r\overline \xi) \prod_{j=1}^n (r_j^2 - R_j^2)^\alpha R_1 \cdots R_n dR d\xi.
\end{align}
\end{lem}

The first part of the following lemma was stated in \cite{Sh2}. This lemma provides direct connection between mixed norm spaces and standard $A^{p}_{\alpha}$ Bergman classes in higher dimensions. A detailed proof of the first part can be found in the cited paper of the second author and the proof of the second part is just a modification of that proof. We again omit details refereing the reader to \cite{Sh2}.

\begin{lem}\label{LeEmb}
Let $0< \max(p,q) \leq s < \infty$ and $\alpha > 0$. Then we have
\begin{equation}\label{EqEmb}
\left( \int_{\mathbb D^n} |f(w)|^s (1-|w|)^{s(\alpha + \frac{1}{p})-2} dV(w) \right)^{1/s} \leq C
\| f \|_{F^{p,q}_\alpha}, \qquad f \in F^{p,q}_\alpha (\mathbb D^n),
\end{equation}
\begin{equation}\label{EqEmb1}
\left( \int_{\mathbb D^n} |f(w)|^s (1-|w|)^{s(\alpha + \frac{1}{p})-2} dV(w) \right)^{1/s} \leq C
\| f \|_{A^{p,q}_\alpha}, \qquad f \in A^{p,q}_\alpha (\mathbb D^n).
\end{equation}
\end{lem}

The following lemma is crucial for all proofs of necessity of multiplier conditions in our main results. It provides explicit estimates of the denominator of the Bergman kernel in various (quasi) norms in polydisc which appear in this note. Again these estimates are known and can be found in \cite{Sh1}, \cite{Sh2}, \cite{Sh3}. Let us note that the Bergman kernel in polydisc is a product of $n$ one dimensional Bergman kernels. This fact often allows one to reduce calculations in several variables to the already classical one dimensional case, see for example \cite{DS} and references therein.

\begin{lem}[\cite{Sh1}]\label{LeSh1}
Let $0<p,q<\infty$ and set
\begin{equation}\label{Eqtestfr}
f_R(z) = \frac{1}{(1-Rz)^{\beta + 1}}, \qquad \beta > -1, \quad R \in I^n, \quad z \in \mathbb D^n.
\end{equation}
Then we have the following (quasi) norm estimates:
\begin{align}\label{EqMZ1}
\| f_R \|_{H^p(\mathbb D^n)} \leq & \frac{C}{(1-R)^{\beta - 1/p + 1}}, \qquad \beta > \frac{1}{p} - 1 \\
\| f_R \|_{A^{p,q}_\alpha (\mathbb D^n)} \leq & \frac{C}{(1-R)^{\beta - \alpha - 1/p + 1}},\qquad \beta >\alpha-1+\frac{1}{p} \\
\| f_R \|_{F^{p,q}_\alpha (\mathbb D^n)} \leq &
\frac{C}{(1-R)^{\beta - \alpha - 1/p + 1}},\qquad \beta
>\alpha-1+\frac{1}{p}
\end{align}
\end{lem}
Note the last estimate of last lemma is true also for two limit spaces $F^{p,\infty,s}$ and $A^{p,\infty,s}$ we defined above

The following lemma follows from Lemma \ref{LeEmb} but we include it here to stress its importance for later proofs.

\begin{lem}[\cite{DS}]\label{LeDS}
If $0 < v \leq 1$, $q>\frac{1}{v}-2$ and $t > 0$ then we have, for every $f \in H(\mathbb D^n)$:
\begin{equation}\label{EqDS}
\int_{\mathbb D^n} |f(w)|^t(1-|w|)^{q} dV(w) \leq C \left( \int_{\mathbb D^n} |f(w)|^{vt} (1-|w|)^{2v - 2+qv} dV(w) \right)^{1/v}.
\end{equation}
\end{lem}
The proof of this lemma shows that if we replace here $(1-|w|)^q$ by $(1-|w|r)^q$, where $r\in I^n$, then the integrand on the right hand side changes to
$$(1-|w|r)^{qv}(1-|w|)^{2v-2}|f(w)|^{vt}$$ and conditions on parameters in this estimate will be $$\frac{1}{2}< v \leq 1$$.
The proof of this last statement is a  modification of the well-known proof of Lemma \ref{LeDS}, therefore we omit easy details.
This remark will be nevertheless used by us in the proof of the main result of this note.

%The last auxiliary result we need is the following lemma.

%\begin{lem}[\cite{DS}]\label{rro}
%Let $\alpha > -1$ and $\lambda > \alpha + 1$. Then
%$$\int_0^1 \frac{(1-R)^\alpha}{(1-Rr)^\lambda} dR \leq C (1-r)^{\alpha + 1 - \lambda}, \qquad 0 \leq r < 1.$$
%\end{lem}

The first assertions we formulate provide conditions which are necessary, but not in general sufficient, in cases when all indexes are different in pairs of mixed norm spaces we consider in this note. By this we mean spaces of coefficient multipliers from $F^{p,q}_{\alpha}$ spaces to $A^{r,s}_\beta$ spaces in the polydisc, and conversely. Even these assertions are more general than those that are present in the literature, see \cite{DS}, a large survey article \cite{Sv1} by V. Shvedenko and more recent work \cite{Bl1}, \cite{BlA}, \cite{BP} by O. Blasco, J. L. Arregui and M. Pavlovic. After that we formulate a theorem containing a sharp result which is one of the main results of this paper. This last result also extends known assertions on coefficient multipliers of Bergman spaces in the unit disk in two directions at the same time: it deals with mixed norm $F^{p,q}_{\alpha}$ spaces and deals with polydisc instead of the disk. Since the proof of our main theorem contains in some sense proofs of all preliminary easy observations which we put below before formulation of main theorem we provide only complete sketches of proofs of these observations.

Let $g\in H(\mathbb D^n)$, $g(z)=\sum_{k \in \mathbb Z_+^n} c_k z^k $. If $c = \{ c_k \}_{k \in \mathbb Z_+^n}$ is a coefficient multiplier from $A^{p,q}_\alpha$ to $F^{t,s}_\beta$ where $t \leq s$, then the following condition holds:
\begin{equation}\label{cond}
\sup_{r\in I^{n}} M_{t}(D^{m}g,r)(1-r)^{\beta+m-\alpha-\frac{1}{p}+1}<\infty,
\end{equation}
where $m$ is any number satisfying $m>\alpha-\beta+\frac{1}{p}-1$.

Moreover the last condition is a necessary condition for the following cases when we consider multipliers from $F^{p,q}_\alpha$ to $A^{t,s}_\beta$ or from $F^{p,q}_\alpha$ to   $F^{t,s}_\beta$ for $t \leq s$.

However, we add additional conditions when we look at multipliers into Lizorkin-Triebel type spaces. These conditions appear due to the first embedding in the following chain
$$F^{t,s}_\beta \hookrightarrow A^{t,s}_\beta \hookrightarrow A^{t,\infty}_{\beta}$$
of embeddings. The first embeddings holds for $t\leq s$ and follows directly from Minkowski inequality. The second embedding is a well known classical fact and it holds for all strictly positive $s$, $t$ and $\beta$, see, for example, \cite{DS}, \cite{OF1}, \cite{OF2}.

In all of the above cases of necessity of condition (\ref{cond}) proof follows directly from Lemma \ref{LeSh1} and the above embeddings in combination with the closed graph theorem which is always used in such situations, see the proof of theorem below as a typical example. We remark in addition that in the case of multipliers into Lizorkin-Triebel spaces we should use both embeddings, in the other case only the second one. We omit  details referring readers to the proof of the first part of
Theorem \ref{Th3IzvVu}.We will see this necessary condition on a multiplier which we just mentioned holds also  when $s\leq t$,moreover we will show that it is sharp in this case.

If $t = \infty$ and $\alpha = \gamma + r$, then the condition
(\ref{cond}) is again necessary when considering multipliers from
$F^{p,q}_\gamma$ into $A^{\infty,\infty}_{r, \beta}$. We again
omit a proof which is a modification of arguments of the proof of
the first part of Theorem \ref{Th3IzvVu}.

\begin{thm}\label{Th3IzvVu}
Let $g \in H(\mathbb D^n)$, $g(z) = \sum_{k \in \mathbb Z_+^n} c_k
z^k$. Assume $t\leq 1$, $\frac{t}{2}< s \leq t < \infty$, $0 <
\max(p,q) \leq s <\infty$, $\beta+\frac{1}{t}<\alpha + \frac{1}{p}
< \frac{2}{t}$,$m \in \mathbb N$ and $m
> \frac{2}{t}-1$.

$1^o$. $\{ c_k \}_{k \in \mathbb Z_+^n} \in M_T(F^{p,q}_\alpha (\mathbb D^n), A^{t,s}_\beta(\mathbb D^n))$ if and only if
\begin{equation}\label{EqTh31}
\sup_{r \in I^n} M_t(D^mg, r)(1-r)^{m+1 - \frac{1}{p} + \beta - \alpha} < \infty.
\end{equation}

$2^o$. $\{ c_k \}_{k \in \mathbb Z_+^n} \in M_T(A^{p,q}_\alpha (\mathbb D^n), A^{t,s}_\beta(\mathbb D^n))$ if and only if
\begin{equation}\label{EqTh32}
\sup_{r \in I^n} M_t(D^mg, r)(1-r)^{m+1 - \frac{1}{p} + \beta - \alpha} < \infty.
\end{equation}
\end{thm}

We remark that the above mentioned embedding $A^{p,q}_\alpha \hookrightarrow F^{p,q}_\alpha$, $p\geq q$ (which follows from Minkowski's inequality) allows us to deduce sufficiency of condition (\ref{EqTh32}) from part $1^o$ of the above theorem under
additional condition $p \geq q$. However, this additional condition can be dropped with help of the second part of Lemma \ref{LeEmb}. Hence we will omit the complete proof of second part providing only the complete sketches of it.

{\it Proof.} $1^o$. Assume $\{ c_k \}_{k \in \mathbb Z_+^n} \in M_T(F^{p,q}_\alpha (\mathbb D^n), A^{t,s}_\beta(\mathbb D^n))$. An application of the closed graph theorem gives $\| M_c f \|_{A^{t,s}_\beta} \leq C \| f \|_{F^{p,q}_\alpha}$. Let $w \in
\mathbb D^n$ and set $f_w(z) = \frac{1}{1-wz}$, $g_w = M_c f_w$. Then we have
$$D^m g_w = D^m M_c f_w = M_c D^m f_w = C M_c \frac{1}{(1-wz)^{m+1}},$$
which, together with the third estimate from Lemma \ref{LeSh1}, gives
\begin{equation}
\| D^m g_w \|_{A^{t,s}_\beta} \leq C \left\| \frac{1}{(1-wz)^{m+1}} \right\|_{F^{p,q}_\alpha}
\leq \frac{C}{(1-|w|)^{m-\alpha - 1/p +1}}.
\end{equation}
Now Lemma \ref{LeEmb}, where $f$, $p$, $q$, $s$ and $\alpha$ are replaced respectively by $D^m g_w$, $t$, $s$, $t$ and $\beta$, combined with the above estimate, gives
\begin{equation}
\left( \int_{\mathbb D^n} |D^m g_w(z)|^t (1-|z|)^{t\beta - 1} dV(z) \right)^{1/t} \leq \frac{C}{(1-|w|)^{m - \alpha - 1/p + 1}}.
\end{equation}
Setting $|w| = r$ and using estimate
\begin{equation}\label{EqMtEst}
M_t(\phi, r)(1-r)^\beta \leq C \left( \int_{\mathbb D^n} |\phi(z)|^t (1-|z|)^{t\beta - 1} dV(z) \right)^{1/t}
\end{equation}
we derive (\ref{EqTh31}).
The above argument provides also a complete proof of necessity in part $2^o$ of the theorem. Indeed, the only change is that
we use the second estimate in Lemma \ref{LeSh1} instead of the third one. Moreover, similar arguments can be used to prove various necessity assertions we stated before this theorem.

Conversely, assume (\ref{EqTh31}) holds and choose $f \in F^{p,q}_\alpha (\mathbb D^n)$. Set $h = M_c f$ and choose
$\eta \in \mathbb T^n$, $r \in I^n$. Then we have, using Lemma \ref{LeCon}:
\begin{align}
h(r^2 \eta) & = \int_{\mathbb T^n} f(r\eta \overline t) g(rt) dt \\
& = \int_{\mathbb T^n} D^m g(R\xi) f(r\eta \overline \xi) \int_0^{r_1} \cdots \int_0^{r_n} \prod_{j=1}^n (r_j^2 - R_j^2)^{m-1}
R_1 \ldots R_n dR d\xi.
\end{align}
The last expression can be transformed by a change of variables $R = r \rho$ and we obtain an estimate
\begin{align}\label{EqTrans}
|h(r^2 \eta)| \leq & C \int_{\mathbb D^n} |(D^m g)(r\rho \xi) f(r\rho\eta \overline \xi)|
\prod_{j=1}^n (1- \rho_j^2)^{m-1} d\rho d\xi, \\
= & C \int_{\mathbb D^n} |(D^m g)(r\rho \xi) \overline f(r\rho\eta \overline \xi)|
\prod_{j=1}^n (1- \rho_j^2)^{m-1} d\rho d\xi.
\end{align}
Now we apply Lemma \ref{LeDS} to an analytic function $D^m g(rz) \overline f(r\eta \overline z)$ and obtain
\begin{equation}
|h(r^2 \eta)|^t \leq C \int_{\mathbb D^n} |(D^m g) (r\rho \xi)|^t |f(r\rho \eta \overline \xi)|^t
\prod_{j=1}^n (1- \rho_j^2)^{t(m+1)-2} d\rho d\xi.
\end{equation}
Next we  integrate over $\eta \in \mathbb T^n$. This yields, taking into account (\ref{EqTh31}):
\begin{align}
\phantom{\leq} & \int_{\mathbb T^n} |h(r^2\eta)|^t  d \eta \\
\leq & C \int_{\mathbb T^n} \int_{I^n} \int_{\mathbb T^n}
|(D^mg)(r\rho\xi)|^t |f(r\rho\eta\overline \xi)|^t (1-\rho^2)^{t(m+1)-2} d\rho d\xi d\eta \\
= & C \int_{\mathbb T^n} \int_{I^n} |(D^mg)(r\rho\xi)|^t M_t^t(f,r\rho)  (1-\rho^2)^{t(m+1)-2}
d\rho d\xi \\
\leq & C \int_{\mathbb T^n} \int_{I^n} M_t^t(f, r\rho) (1-r\rho)^{-t(m+1 - \alpha + \beta - 1/p)}  (1-\rho^2)^{t(m+1)-2}
d\rho d\xi \\
\leq & C \int_{\mathbb T^n} \int_{I^n} M_t^t(f, r\rho) (1-r\rho)^{t(\alpha - \beta + 1/p) - 2} d\rho d\xi, \quad t(m+1) > 2 \\
\leq & C \int_{\mathbb D^n} |f(rw)|^t (1-r|w|)^{t(\alpha - \beta + 1/p) - 2} dV(w).
\end{align}

Now we apply Lemma \ref{LeDS} with $v = s/t$ and with $(1-r|w|)^{q}$ instead of $(1-|w|)^{q}$ where
$q=t(\frac{1}{p}-\beta+\alpha)-2$, see comments after Lemma \ref{LeDS}. This gives
\begin{equation}
M_t^s(h, r^2) \leq C \int_{\mathbb D^n} |f(rw)|^s (1-r|w|)^{s(\alpha - \beta +1/p) - 2s/t} (1-|w|)^{2s/t - 2} dV(w).
%\\ \leq & C \int_{\mathbb D^n} |f(rw)|^s (1-r|w|)^{s(\alpha - \beta +1/p)} (1-|w|)^{- 2} dV(w)
\end{equation}
Next we integrate both sides over $r \in I^{n}$. This gives
\begin{align*}
\| h \|_{A^{t,s}_\beta}^s & = \int_{I^n} M_t^s(h, R)(1-R)^{\beta s - 1} dR \\
& \leq C \int_{I^n} (1-R)^{\beta s -1} \int_{\mathbb D^n} |f(Rw)|^s
\frac{(1-|w|)^{2s/t - 2}}{(1-R|w|)^{-s(\alpha - \beta + 1/p)+2s/t}}dV(w) dR \\
& \leq C \int_{I^n} (1-R)^{\beta s-1} \int_{I^n} \int_{\mathbb T^n} |f(rR\xi)|^s
\frac{(1-|w|)^{2s/t - 2}}{(1-R|w|)^{-s(\alpha - \beta + 1/p)+2s/t}} dr dR d\xi \\
& = C \int_{I^n} \int_{I^n} M_s^s(f, rR) (1-R)^{\beta s-1}
\frac{(1-|w|)^{2s/t - 2}}{(1-R|w|)^{-s(\alpha - \beta + 1/p)+2s/t}}dr dR\\
& \leq C \int_{I^n} (1-r)^{2s/t-2} M_s^s(f, r) \int_{I^n}
\frac{(1-R)^{\beta s -1}}{(1-rR)^{s(\beta - \alpha - 1/p)+2s/t}} dR dr\\
& \leq C \int_{I^n} M_s^s(f, r) (1-r)^{s(\alpha + 1/p) - 2} dr.
\end{align*}
At the last the assumption $t(\alpha + 1/p)<2$ allowed us to use the following well known estimate:
$$\int_0^1 \frac{(1-R)^\alpha}{(1-Rr)^\lambda} dR \leq C (1-r)^{\alpha + 1 - \lambda}, \qquad 0 \leq r < 1,$$
valid for $\alpha > -1$ and $\lambda > \alpha + 1$. Finally, Lemma \ref{LeEmb}, specifically embedding (\ref{EqEmb}), allows us to conclude that $\| h \|_{A^{t,s}_\beta} \leq C \| f \|_{F^{p,q}_\alpha}$.

The proof of sufficiency in the $A^{p,q}_\alpha$ case goes along the same lines, the only difference being use of embedding
(\ref{EqEmb1}) instead of (\ref{EqEmb}) at the last step. $\Box$

%We remark using at last step above the second part of Lemma \ref{LeEmb} we easily see this proof of main theorem is at the %same time  a proof of a sharp result on multipliers into  $A^{t,s}_\beta$ spaces from $A^{p,q}_\alpha$ function spaces. This %is
The following theorem is parallel to the previous one, it gives characterization of multipliers into $F^{t, s}_\beta$ spaces.
The proof goes along similar lines as the proof of the above theorem, therefore we omit it.
Note here the necessity follows directly from the same arguments as in  previous theorem ,but it is based on first estimate
in second lemma.The sufficiency part also repeats the proof of previous theorem till the end, at last step we use the mentioned embedding between $F^{t,s}_\alpha$ and $A^{t,s}_\alpha$

\begin{thm}\label{Th3IzvVuF}
Let $g \in H(\mathbb D^n)$, $g(z) = \sum_{k \in \mathbb Z_+^n} c_k
z^k$. Assume $t\leq1$,$\frac{t}{2}< s \leq t < \infty$, , $0 <
\max(p,q) \leq s <\infty$, $\beta+\frac{1}{t}<\alpha + \frac{1}{p}
< \frac{2}{t}$, $m \in \mathbb N$ ,$m
> \frac{2}{t}-1$.

$1^o$. $\{ c_k \}_{k \in \mathbb Z_+^n} \in M_T(F^{p,q}_\alpha (\mathbb D^n), F^{t,s}_\beta(\mathbb D^n))$ if and only if
\begin{equation}\label{EqTh33}
\sup_{r \in I^n} M_t(D^mg, r)(1-r)^{m+1 - \frac{1}{p} + \beta - \alpha} < \infty.
\end{equation}

$2^o$. $\{ c_k \}_{k \in \mathbb Z_+^n} \in M_T(A^{p,q}_\alpha (\mathbb D^n), F^{t,s}_\beta(\mathbb D^n))$ if and only if
\begin{equation}\label{EqTh34}
\sup_{r \in I^n} M_t(D^mg, r)(1-r)^{m+1 - \frac{1}{p} + \beta - \alpha} < \infty.
\end{equation}
\end{thm}

We remark partially results of this note were announced before in \cite{Sh2} without proofs.

We conclude this paper with remarks on limit case spaces $F^{p, \infty, \alpha}$, which are defined analogously to the
$F^{p,q}_\alpha$ spaces, the only difference being replacement of inner (quasi) norm $L^q$ by $L^\infty$ norm. We also constantly use in arguments an obvious fact that weighted Hardy spaces $A^{p,\infty}_\alpha$ contain $F^{p,\infty,\alpha}$ classes in the polydisc.

Namely, we look for estimates on the rate of growth of a function $g$ which represents a coefficient multiplier into the
above described space. Again proofs are parallel to the one presented in the proof of Theorem \ref{Th3IzvVu}.

We again assume that $g = \sum_{k\in \mathbb Z_+^n} c_k z^k$ is analytic in polydisc, by an estimate of the rate of growth of $g$ we mean the following:
\begin{equation}\label{EqRgr}
\sup_{r\in I^n}M_p(D^{m}g,r)(1-r)^{\tau}<\infty.
\end{equation}
The condition (\ref{EqRgr}) is often a necessary condition for $\{ c_k \}_{k \in \mathbb Z_+^n}$ to be a multiplier and in
that case $\tau$ depends on parameters involved and can be explicitly specified for each pair of spaces separately using Lemma \ref{LeSh1}. A typical example is the following
\begin{prop}
Assume $0< v \leq p$ and $m > \frac{1}{p} -1-s$.
If a sequence $\{ c_k \}_{k \in \mathbb Z_+^n}$ is a multiplier from $H^v$ to $F^{p,\infty,s}$, then the function $g =
\sum_{k \in \mathbb Z_+^n} c_k z^k$ satisfies condition (\ref{EqRgr}) with $\tau = m + 1 + s - 1/p$.
\end{prop}
The same  condition  on $g$ function, but for different $\tau$ is necessary for $g$ to be a multiplier from mixed norm spaces we studied in this paper to $F^{p,\infty, s}$ or from $F^{q,\infty,s}$ to $F^{p,\infty,v}$. In the case of multipliers from $F^{p,\infty,v}$ into $A^{\infty,\infty}_{s,t}$ classes the necessary condition is the same (for different $\tau$) but with $\infty$ instead of $p$. We leave simple proofs of these facts based on detailed outlines provided above to interested readers.

\end{document}